\documentclass[11pt,a4paper]{article}
\usepackage{pdfsync} %

\usepackage{amsfonts,amsmath,amsthm}
\newtheorem{theorem}{Theorem}
\numberwithin{theorem}{section}
\newtheorem{lemma}[theorem]{Lemma}
\newtheorem{corollary}[theorem]{Corollary}
\newtheorem{example}[theorem]{Example}
\usepackage{amssymb}
\usepackage{hyperref}
\usepackage{graphicx}
\usepackage{wrapfig}
\usepackage{subfig}
\usepackage{epsfig}
\usepackage{float,mathrsfs}
\usepackage{enumitem}
\usepackage{bbm}

\usepackage{algpseudocode}
\usepackage[font=scriptsize]{caption}
\usepackage{chngcntr} 
\counterwithout{equation}{section}
\counterwithout{figure}{section}
\usepackage{color}
\definecolor{refkey}{rgb}{0.9451,0.2706,0.4941}\definecolor{labelkey}{rgb}{0.9451,0.2706,0.4941}

\definecolor{darkred}{RGB}{139,0,0}
\definecolor{darkgreen}{RGB}{0,100,0}
\definecolor{darkmagenta}{RGB}{139,0,139}

\newcommand{\bsx}{{\boldsymbol{x}}}
\newcommand{\bsm}{{\boldsymbol{m}}}
\newcommand{\bsmu}{{\boldsymbol{\mu}}}
\newcommand{\bsy}{{\boldsymbol{y}}}
\newcommand{\bsalpha}{{\boldsymbol{\alpha}}}

\newcommand{\bst}{{\boldsymbol{t}}}

\newcommand{\bsgamma}{{\boldsymbol{\gamma}}}

\newcommand{\bbN}{{\mathbb{N}}}
\newcommand{\bbR}{{\mathbb{R}}}

\newcommand{\bsnu}{{\boldsymbol{\nu}}}

\newcommand{\bsb}{{\boldsymbol{b}}}

\numberwithin{equation}{section}

\title{Generalized dimension truncation error analysis\\for high-dimensional numerical integration:\\lognormal setting and beyond}

\author{Philipp A.~Guth\footnotemark[2]\and Vesa Kaarnioja\footnotemark[3]}

\renewcommand{\thefootnote}{\fnsymbol{footnote}}

\begin{document}
\maketitle

\begin{abstract}
Partial differential equations (PDEs) with uncertain or random inputs have been considered in many studies of uncertainty quantification. In forward uncertainty quantification, one is interested in analyzing the stochastic response of the PDE subject to input uncertainty, which usually involves solving high-dimensional integrals of the PDE output over a sequence of stochastic variables. In practical computations, one typically needs to discretize the problem in several ways: approximating an infinite-dimensional input random field with a finite-dimensional random field, spatial discretization of the PDE using, e.g., finite elements, and approximating high-dimensional integrals using cubatures such as quasi-Monte Carlo methods. In this paper, we focus on the error resulting from dimension truncation of an input random field. We show how Taylor series can be used to derive theoretical dimension truncation rates for a wide class of problems and we provide a simple checklist of conditions that a parametric mathematical model needs to satisfy in order for our dimension truncation error bound to hold. Some of the novel features of our approach include that our results are applicable to non-affine parametric operator equations, dimensionally-truncated conforming finite element discretized solutions of parametric PDEs, and even compositions of PDE solutions with smooth nonlinear quantities of interest. As a specific application of our method, we derive an improved dimension truncation error bound for elliptic PDEs with lognormally parameterized diffusion coefficients. Numerical examples support our theoretical findings.
\end{abstract}

\footnotetext[2]{Johann Radon Institute for Computational and Applied Mathematics, Austrian Academy of Sciences, Altenbergerstra{\ss}e 69, A-4040 Linz, Austria, {\tt philipp.guth@ricam.oeaw.ac.at}}
\footnotetext[3]{Department of Mathematics and Computer Science, Free University of Berlin, Arnimallee 6, 14195 Berlin, Germany, {\tt vesa.kaarnioja@fu-berlin.de}}

\renewcommand{\thefootnote}{\arabic{footnote}}

\section{Introduction}
It is common in the field of uncertainty quantification to consider mathematical models where uncertain inputs are parameterized by infinite sequences of random variables. A typical problem considered in the literature is the problem of finding $u:D\times U\to \bbR$ such that
\begin{align}
\begin{split}
-\nabla\cdot(a(\bsx,\bsy)\nabla u(\bsx,\bsy))&=f(\bsx), \quad\bsx\in D,~\bsy\in U,\\
u(\bsx,\bsy)&=0, \quad\quad\bsx\in\partial D,~\bsy\in U,
\end{split}\label{eq:pdeproblem0}
\end{align}
in some bounded Lipschitz domain $D\subset\bbR^d$, $d\in \{1,2,3\}$, where $a\!:D\times U\to \mathbb R$ is a parameterized diffusion coefficient and $f\!:D\to\mathbb R$ is a given source term. The parameter set $U$ is assumed to be a nonempty subset of $\mathbb R^{\mathbb N}$. More generally, consider an abstract mathematical model $M\!:X\times U\to Y$ such that 
\begin{align*}
M(g(\bsy),\bsy)=0,
\end{align*}
where $X$ and $Y$ are separable Banach spaces and $U$ is defined as above. A natural quantity to investigate is the expected value
\begin{align}
I(g):=\int_U g(\bsy)\,\bsmu({\rm d}\bsy),\label{eq:infobj}
\end{align}
where $\bsmu$ is a probability measure over $U$. In many applications, $\bsmu$ is either chosen as the uniform probability measure over $U=[-1,1]^{\mathbb N}$ or a Gaussian probability measure over $U=\mathbb R^{\mathbb N}$. The model $M$ can for instance represent the PDE problem~\eqref{eq:pdeproblem0}, in which case one could consider either $g(\bsy)=u(\cdot,\bsy)$ (in which case the integral~\eqref{eq:infobj} is interpreted as a Bochner integral) or $g(\bsy)=G(u(\cdot,\bsy))$, where $G$ is either a linear or nonlinear functional (quantity of interest). However, we stress that the function $\bsy\mapsto g(\bsy)$ could equally well represent the output of any kind of parametric model such as a parabolic PDE~\cite{guth22}, a hyperbolic PDE~\cite{hoangschwab}, or even the eigenpairs of a spectral eigenvalue problem~\cite{GGKSS2019} with uncertain input parameters, as well as the composition of these outputs with suitable (possibly nonlinear) quantities of interest.

For the numerical treatment of~\eqref{eq:infobj}, a natural first step would be to consider a dimensionally-truncated model $M_s\!:X\times U_s\to Y$ such that
\begin{align}
M_s(g_s(\bsy_{\leq s}),\bsy_{\leq s})=0,\label{eq:modelproblemM2}
\end{align}
where $\varnothing\neq U_s\subseteq \mathbb R^s$ and $g_s(\bsy_{\leq s})\in X$ for all $\bsy_{\leq s}\in U_s$. The corresponding expected value in this case is then given by
$$
I_s(g_s):=\int_{U_s}g_{s}(\bsy_{\leq s})\,\bsmu_{\leq s}({\rm d}\bsy_{\leq s}),
$$
where $\bsmu_{\leq s}$ denotes an appropriate probability measure on $U_s$. By considering $I_s(g_s)$ instead of $I(g)$, we have introduced a \emph{dimension truncation error}
$$
\|I(g)-I_s(g_s)\|_X.
$$

In many practical problems involving partial differential equations (PDEs), there are also other sources of errors: for example, the integral operator $I_s$ may need to be approximated by a cubature rule $Q_{s,n}$ with $n$ nodes and, in practice, we may only have access to an approximation $g_{s,h}$ of the solution to~\eqref{eq:modelproblemM2} in some finite-dimensional subspace $X_h$ of $X$ (for example, a conforming finite element subspace). The overall error can typically be estimated via an error decomposition of the form
$$
\|I(g)-Q_{s,n}(g_{s,h})\|_X\leq \|I(g)-I_s(g_s)\|_X+\|I_s(g_s-g_{s,h})\|_X+\|I_s(g_{s,h})-Q_{s,n}(g_{s,h})\|_X,
$$
where the last two terms correspond to \emph{spatial discretization error} and \emph{cubature error}, respectively.

In this work, we focus on the analysis of the dimension truncation error for integral quantities. The contribution of dimension truncation error is independent of the numerical scheme chosen for the cubature operator or spatial discretization, which allows us to approach the problem from a general vantage point---for example, we do not need to restrict our analysis to a specific numerical cubature method, spatial discretization scheme or even a specific mathematical model problem. Instead, we proceed to derive general conditions under which it is possible to derive explicit rates for the dimension truncation error.

Dimension truncation error rates are typically studied in conjunction with elliptic PDEs with random inputs. Kuo \emph{et al}.~\cite{kss12} derived a dimension truncation rate for a class of elliptic PDEs under an affine parameterization of the uncertain diffusion coefficient. This result was later improved by Gantner~\cite{gantner}, who studied dimension truncation in the context of affine parametric operator equations. Dimension truncation has also been studied for coupled PDE systems arising in optimal control problems under uncertainty~\cite{guth21} as well as in the context of the so-called ``periodic model'' of uncertainty quantification for both numerical integration~\cite{kks20} and kernel interpolation~\cite{kkkns22}. A common feature in these works is the use of Neumann series, which is a suitable tool for dimension truncation analysis provided that the uncertain parameters enter the PDE operator affinely. A drawback of this approach is that when the dependence of the PDE operator on the random variables becomes sufficiently nonlinear, the Neumann approach may produce only suboptimal dimension truncation rates: this is the case for lognormally parameterized diffusion coefficients for elliptic PDEs~\cite{log} or when the quantity of interest is a nonlinear functional of the PDE response~\cite{dglgs19,hks21}.

In contrast to using Neumann series to obtain dimension truncation error rates, the use of Taylor series was considered by Gilbert \emph{et al}.~\cite{GGKSS2019} in the context of a spectral eigenvalue problem for an elliptic PDE with a random coefficient. The use of Taylor series makes it possible to exploit the \emph{parametric regularity} of the model problem in order to derive dimension truncation rates. %
 A similar Taylor series approach---motivated by the paper~\cite{GGKSS2019}---was used to derive a dimension truncation rate for a smooth nonlinear quantity of interest subject to an affine parametric parabolic PDE in~\cite{guth22}.

The available literature provides some numerical evidence concerning dimension truncation rates for nonlinear parameterizations of PDE problems. The PhD thesis~\cite{gantnerphd} provides some numerical experiments suggesting that the dimension truncation error rates for certain non-affine parametric PDE problems are actually significantly better than the theoretical bounds derived using Neumann series. Both the numerical and theoretical results in~\cite{guth22} indicate that the use of a smooth nonlinear quantity of interest does not deteriorate the dimension truncation rate for an affine parametric PDE problem. Furthermore, it is known in the context of lognormal parameterizations for diffusion coefficients of parametric elliptic PDEs that the use of special Mat\'ern covariances can yield even exponentially convergent dimension truncation errors (cf.~\cite[Section~7.2]{matern1} and~\cite{matern2,schwabtodor}).

The contributions of our paper are the following. We derive a dimension truncation error bound for generic, parameterized elements of a separable Banach space under moderate assumptions on parametric regularity and integrability with respect to $\beta$-Gaussian probability measures. Especially, the special case $\beta=2$ corresponds to the Gaussian probability measure and the limiting case $\beta\to\infty$ corresponds to the uniform probability measure. Our approach to analyzing the dimension truncation error has several key advantages:
\begin{itemize}
    \item The results can be applied to parametric PDEs by simply checking that the conditions of the dimension truncation theorem hold in their respective Sobolev spaces. Especially, the dimension truncation rate we obtain as a corollary for elliptic PDEs with lognormally parameterized diffusion coefficients improves the existing results~\cite{log}.
    \item Since the results are stated for general separable Banach spaces, we may also apply our results for conforming finite element discretized solutions of parametric PDEs---a result which appears to have received little attention in the existing literature save for~\cite[Theorem~8]{log}. This is made possible by the fact that conforming finite element solutions inherit the parametric regularity of the PDE problem as opposed to, e.g., non-conforming finite element methods or other discretization schemes.
    \item Finally, since our analysis is not restricted to only parametric PDE problems, the results are also valid for, e.g., the composition of the solution of a parametric PDE with a smooth nonlinear quantity of interest, provided that the composition satisfies the hypotheses of our dimension truncation theorem. We illustrate this by giving a specific example of a certain nonlinear quantity of interest with known parametric regularity bounds.
    \end{itemize}

This document is structured as follows. Subsection~\ref{sec:notations} introduces the multi-index notation used throughout the paper and Stechkin's lemma is stated for the reader's convenience. The problem setting is introduced in Section~\ref{sec:problem}, containing the fundamental assumptions needed for the ensuing dimension truncation analysis. We briefly discuss the concept of infinite-dimensional integration in Section~\ref{sec:infint}. Section~\ref{sec:dimtrunc} contains the main dimension truncation theorems developed in this work, and the consequences for parametric elliptic PDE problems are outlined in Section~\ref{sec:pdeappl}. Numerical examples assessing the sharpness of our theoretical results are presented in Section~\ref{sec:numex}. The paper ends with some conclusions.

\subsection{Notations and preliminaries}\label{sec:notations}
Throughout this manuscript, boldfaced letters are used to denote multi-indices while the subscript notation $m_j$ is used to refer to the $j^{\rm th}$ component of multi-index $\boldsymbol m$. We denote the set of all finitely supported multi-indices by
$$
\mathcal F := \{\boldsymbol m\in\mathbb N_0^{\mathbb N}:|{\rm supp}(\boldsymbol m)|<\infty\},
$$
where the support of a multi-index is defined as ${\rm supp}(\boldsymbol m):=\{j\in\mathbb N:m_j\neq 0\}$. Moreover, the modulus of a multi-index is defined as
$$
|\boldsymbol m|:=\sum_{j\geq 1} m_j.
$$
Furthermore, for any sequence $\boldsymbol x:=(x_j)_{j=1}^\infty$ of real numbers and $\boldsymbol m\in\mathcal F$, we define the notation
$$
\boldsymbol x^{\boldsymbol m}:=\prod_{j\geq 1} x_j^{m_j},
$$
where we use the convention $0^0:=1$.

The following lemma is commonly used in the analysis of best $N$-term approximations (cf., e.g.,~\cite{devore}), and it will be highly useful in our treatment of the dimension truncation error.
\begin{lemma}[Stechkin's lemma]\label{lemma:stechkin}
Let $0<p\leq q<\infty$ and let $(a_k)_{k\geq 1}\in\ell^p(\mathbb N)$ be a sequence of real numbers ordered such that $|a_1|\geq |a_2|\geq \cdots$. Then
$$
\bigg(\sum_{k>N}^\infty |a_k|^q\bigg)^{\frac{1}{q}}\leq N^{-\frac{1}{p}+\frac{1}{q}}\bigg(\sum_{k\geq 1} |a_k|^p\bigg)^{\frac{1}{p}}.
$$
\end{lemma}
\proof For an elementary proof of this result, see, e.g., \cite[Lemma 3.3]{KressnerTobler}.\quad\endproof

\section{Problem setting}\label{sec:problem}

Let $\bsalpha:=(\alpha_j)_{j\geq 1}\in\ell^1(\mathbb N)$ be an arbitrary sequence such that $\alpha_j\in[0,\infty)$ for all $j\in\mathbb N$. We define the set
\begin{align*}
    U_{\bsalpha} := \bigg\{ \bsy \in \bbR^{\bbN}: \sum_{j\geq 1} \alpha_j |y_j| <\infty\bigg\}.
\end{align*}
Let $X$ be a separable Banach space, let $g:U_{\bsalpha}\to X$, and 
let us define $g_s(\bsy):=g(\bsy_{\leq s},\mathbf 0):=g(y_1,\ldots,y_s,0,0,\ldots)$. We consider 
\begin{align}
    \bsmu_\beta(\mathrm d\bsy) := \bigotimes_{j\geq1} \mathcal{N}_\beta(0,1),\label{eq:gengauss}
\end{align}
where $\mathcal{N}_\beta(0,1)$ denotes the univariate $\beta$-Gaussian distribution with density 
\begin{align*}
    \varphi_{\beta}(y) := \frac{1}{2\beta^{\frac{1}{\beta}}\Gamma(1+\frac{1}{\beta})} {\rm e}^{-\frac{|y|^\beta}{\beta}}, \quad y\in \bbR,
\end{align*}
where we restrict to the case $\beta \geq1$. Importantly, in the case $\beta=2$ the probability measure~\eqref{eq:gengauss} is Gaussian and in the case $\beta=1$ it corresponds to the Laplace distribution. Formally, the case $\beta=\infty$ corresponds to the uniform probability measure on $[-1,1]^{\mathbb N}$, which we denote by
$$
\bsgamma({\rm d}\bsy):=\bigotimes_{j\geq 1}\frac{{\rm d}y}{2}.
$$
To this end, we will consider the \emph{dimension truncation error}
\begin{align*}
    \bigg\|\int_{\bbR^{\bbN}} (g(\bsy)-g_s(\bsy))\,\bsmu_\beta({\mathrm d}\bsy)\bigg\|_X
\end{align*}
subject to $\beta$-Gaussian probability measures and the uniform probability measure, equipped with their respective sets of assumptions.

Integration problems involving measures of this type occur in Bayesian inverse problems governed by PDEs endowed with Besov priors, see, e.g., \cite{DashtiHarrisStuart12,hks21,schwabstein23}.

\paragraph{$\beta$-Gaussian probability measures.} In the $\beta$-Gaussian setting, we will work under the following assumptions:
\begin{enumerate}[label=(A\arabic*)]
\item It holds for a.e.~$\bsy\in U_{\bsalpha}$ that\label{assump1}
$$
\|g(\bsy)-g_s(\bsy)\|_X \to 0 \quad \text{as }s\to \infty.
$$
\item Let $(\Theta_k)_{k\geq 0}$ be a sequence of nonnegative numbers, let $\bsb := (b_j)_{j\geq 1} \in \ell^p(\bbN)$ for some $p \in (0,1)$, and let $b_1 \geq b_2 \geq \cdots \geq 0$. We assume that the integrand $g$ is continuously differentiable up to order $k+1$, with
    \begin{align*}
        \|\partial^{\bsnu}g(\bsy)\|_X \leq \Theta_{|\bsnu|} \bsb^{\bsnu} \prod_{j \geq 1} {\rm e}^{\alpha_j|y_j|}<\infty
    \end{align*}
    for all $\bsy \in U_{\bsalpha}$ and all $\bsnu \in \mathcal{F}_k := \{\bsnu \in \bbN_0^\bbN: |\bsnu|\leq k+1\}$, where $k := \lceil\frac{1}{1-p}\rceil$. In the case $\beta = 1$, we assume in addition that $\alpha_j < 1$ for all $j \in \bbN$.\label{assump2}
\end{enumerate}

\paragraph{Uniform probability measure.} In this setting, we suppose that $g(\bsy)\in X$ for each $\bsy\in[-1,1]^{\mathbb N}$ and we will work under the following assumptions:
\begin{enumerate}[label=(A\arabic*')]
\item It holds for a.e.~$\bsy\in [-1,1]^{\mathbb N}$ that\label{assump1pr}
$$
\|g(\bsy)-g_s(\bsy)\|_X \to 0 \quad \text{as }s\to \infty.
$$
\item Let $(\Theta_k)_{k\geq 0}$ be a sequence of nonnegative numbers, let $\bsb := (b_j)_{j\geq 1} \in \ell^p(\bbN)$ for some $p \in (0,1)$, and let $b_1 \geq b_2 \geq \cdots \geq 0$. We assume that the integrand $g$ is continuously differentiable up to order $k+1$, with
    \begin{align*}
        \|\partial^{\bsnu}g(\bsy)\|_X \leq \Theta_{|\bsnu|} \bsb^{\bsnu}
    \end{align*}
    for all $\bsy \in [-1,1]^{\mathbb N}$ and all $\bsnu \in \mathcal{F}_k := \{\bsnu \in \bbN_0^\bbN: |\bsnu|\leq k+1\}$, where $k := \lceil\frac{1}{1-p}\rceil$.\label{assump2pr}
\end{enumerate}

\emph{Remark.}~(i) One can think of the relationship between $\boldsymbol\alpha$, $g$, and $\boldsymbol b$ as follows: fixing any $\boldsymbol\alpha\in\ell^1(\mathbb N)$ determines the domain $U_{\boldsymbol\alpha}$ of $g$. Once $g$ is fixed, the assumption~\ref{assump2} is required to hold for some $\boldsymbol b\in\ell^p(\mathbb N)$.

(ii) Bounds on the partial derivatives as in (A2) and (A2') typically appear in the context of parameterized PDE models such as those discussed in Section~\ref{sec:pdeappl}.

(iii) Certain holomorphic functions admit bounds of this form, see, e.g., \cite[Proposition 2.3]{SchwabZech}. From~\cite[Proposition~1.2.34]{ZechPhD} we know that compositions of holomophic functions are again holomorphic.

(iv) The summability of the sequence $\bsb$ in \ref{assump2} and \ref{assump2pr} directly relates to the smoothness of the problem. For instance, in case of a holomorphic parametric map $g$, $\bsb$ controls the radii of the domains of analytic continuation and the summability exponent $p$ formalizes a notion for the decay of the $\bsnu$-th partial derivatives of $g$, see~\cite{SchwabZech} and the references therein. The smaller the value of $p$, the faster the decay rate, which will be reflected in the dimension truncation error rates in Theorems~\ref{thm:main} and~\ref{thm:main2}.

\section{Infinite-dimensional integration}\label{sec:infint}

Let $g(\bsy)\in X$ for each $\bsy\in U_{\bsalpha}$. If \ref{assump2} holds, then we infer that $\bsy \mapsto G(g(\bsy))$ for all $G\in X'$ is continuous as a composition of continuous mappings. Hence $\bsy \mapsto G(g(\bsy))$ is measurable for all $G \in X'$, i.e., $\bsy \mapsto g(\bsy)$ is weakly measurable. Since $X$ is assumed to be a separable Banach space, by Pettis' theorem (cf., e.g., \cite[Chapter 4]{Yosida}) we obtain that $\bsy \mapsto g(\bsy)$ is strongly measurable. The $\bsmu_{\beta}$-integrability of the upper bound in \ref{assump2} is proved in \cite[Proposition 3.2]{hks21} for $\beta \in [1,2]$ and can be proved \emph{mutatis mutandis} for $\beta>2$. Thus we conclude from Bochner's theorem (cf., e.g., \cite[Chapter 5]{Yosida}) and \ref{assump2} that $g$ is $\bsmu_\beta$-integrable over $U_{\bsalpha}$. Bochner's theorem can also be used to ensure that a function $g(\bsy)\in X$, $\bsy\in [-1,1]^{\mathbb N}$, is $\bsgamma$-integrable provided that \ref{assump2pr} holds.

The following lemma has been adapted from~\cite[Lemma 2.28]{schwab_gittelson_2011} to our setting.
\begin{lemma}
It holds that $U_{\bsalpha} \in \mathcal{B}(\bbR^\bbN)$, where $\mathcal{B}(\mathbb R^{\mathbb N})$ denotes the Borel $\sigma$-algebra generated by the probability measure $\bsmu_{\beta}$ and the Borel cylinders in $\mathbb R^{\mathbb N}$. Moreover, $\bsmu_\beta(U_{\bsalpha}) = 1$.
\end{lemma}

\proof
    The first statement follows from 
    \begin{align*}
        U_{\bsalpha} = \bigcup_{N\geq 1} \bigcap_{M\geq 1} \bigg\{\bsy \in \bbR^\bbN: \sum_{1\leq j\leq M} \alpha_j |y_j| \leq N\bigg\}.
    \end{align*}
    By the monotone convergence theorem, we deduce that
    \begin{align*}
        \int_{\bbR^\bbN} \sum_{j\geq 1} \alpha_j |y_j|\, \bsmu_{\beta}(\mathrm d\bsy) = \sum_{j\geq 1} \alpha_j \int_{\bbR^\bbN} |y_j|\, \bsmu_{\beta}(\mathrm d\bsy) = \frac{\Gamma(\frac{2}{\beta})}{\beta^{1-\frac{1}{\beta}}\Gamma(1+\frac{1}{\beta})} \sum_{j\geq1} \alpha_j < \infty, 
    \end{align*}
    for all $\beta >0$, where we used~\cite[formula 3.326.2]{gradshteynryzhik}.\quad\endproof

From the above lemma we conclude that we can restrict to $\bsy \in U_{\bsalpha}$ since
\begin{align*}
    \int_{\bbR^\bbN} g(\bsy) \,\bsmu_\beta(\mathrm d\bsy) = \int_{U_\bsalpha} g(\bsy) \,\bsmu_\beta(\mathrm d\bsy),
\end{align*}
for any $g$ satisfying \ref{assump2}. Thus in the $\beta$-Gaussian setting, the domain of integration $\mathbb R^{\mathbb N}$ is interchangeable with $U_{\bsalpha}$.

\begin{lemma}[{\cite[Theorem 1]{infiniteintegration} and \cite[Section 26]{Halmos}}]\label{lemma:infinite}
Let $U\subseteq \bbR^\bbN$, let $\bsmu$ be a probability measure, and let $F:U\to \bbR$ be $\bsmu$-integrable, which satisfy
\begin{align*}
    \lim_{s\to\infty} F(\bsy_{\leq s},\boldsymbol{0}) = F(\bsy) \quad\text{for}~a.e.~\bsy
\intertext{and}
    |F(\bsy_{\leq s},\boldsymbol{0})| \leq |h(\bsy)|\quad\text{for}~a.e.~\bsy 
\end{align*}
for some $\bsmu$-integrable $h$.
From Lebesgue's dominated convergence theorem, we infer the following results.
\begin{itemize}
    \item[\rm (i)] If $U=\bbR^\bbN$ and $\bsmu = \bsmu_{\beta}$, then 
\begin{align*}
    \lim_{s\to \infty} \int_{\bbR^s} F(\bsy_{\leq s},\boldsymbol{0})\,\bsmu_\beta(\mathrm d\bsy_{\leq s}) = \int_{\bbR^\bbN} F(\bsy)\,\bsmu_\beta(\mathrm d\bsy).
\end{align*}
\item[\rm (ii)] If $U=[-1,1]^\bbN$ and $\bsmu = \bsgamma$, then 
\begin{align*}
    \lim_{s\to \infty} \int_{[-1,1]^s} F(\bsy_{\leq s},\boldsymbol{0})\,\bsgamma(\mathrm d\bsy_{\leq s}) = \int_{[-1,1]^\bbN} F(\bsy)\,\bsgamma(\mathrm d\bsy).
\end{align*}
\end{itemize}
\end{lemma}

\emph{Remark.}
In the special case where $F_G(\bsy) := \langle G,g(\bsy) \rangle_{X',X}$ for some $G$ in the topological dual space of $X$, it holds that
\begin{align*}
    F_G(\bsy) - F_G(\bsy_{\leq s},\boldsymbol{0}) = \langle G,g(\bsy) - g(\bsy_{\leq s},\boldsymbol{0})\rangle_{X',X} \leq \|G\|_{X'} \|g(\bsy) - g(\bsy_{\leq s},\boldsymbol{0})\|_X,
\end{align*}
and
\begin{align*}
    |F_G(\bsy_{\leq s},\boldsymbol{0})| \leq \|G\|_{X'} \|g(\bsy_{\leq s},\boldsymbol{0})\|_X,
\end{align*}
which can be bounded by taking $\bsnu = \boldsymbol{0}$ in~\ref{assump2} or~\ref{assump2pr}, respectively. Thus, Lemma~\ref{lemma:infinite} holds for~$F=F_G$ due to~\ref{assump1} in the $\beta$-Gaussian setting and due to~\ref{assump1pr} in the uniform setting.

\section{Dimension truncation error}\label{sec:dimtrunc} The main dimension truncation result is given below.
\begin{theorem}\label{thm:main} Let $g(\bsy)\in X$ for all $\bsy\in U_{\bsalpha}$. Suppose that assumptions {\rm\ref{assump1}} and {\rm\ref{assump2}} hold. Then
$$
\bigg\|\int_{\mathbb R^{\mathbb N}}(g(\bsy)-g_s(\bsy))\bsmu_\beta({\mathrm d}\bsy)\bigg\|_X\leq Cs^{-\frac{2}{p}+1},
$$
where the constant $C>0$ is independent of the dimension $s$.

Let $G\in X'$ be arbitrary. Then
$$
\bigg|\int_{\mathbb R^{\mathbb N}}G(g(\bsy)-g_s(\bsy))\bsmu_\beta({\mathrm d}\bsy)\bigg|\leq C\|G\|_{X'}s^{-\frac{2}{p}+1},
$$
where the constant $C>0$ is above.
\end{theorem}
\proof Let $s^*$ be the smallest integer such that $\sum_{j>s^*}b_j\leq \frac{1}{2}$. Clearly,
$$
\bigg\|\int_{\mathbb R^{\mathbb N}}(g(\bsy)-g_s(\bsy))\,\bsmu_\beta({\rm d}\bsy)\bigg\|_X\leq 2\Theta_0 \prod_{j\geq 1}\int_{\mathbb R}{\rm e}^{\alpha_j|y_j|}\varphi_\beta(y_j)\,{\rm d}y_j=:C_{\rm init}<\infty,
$$
where we used Lemma~\ref{lemma:infinite} and Fubini's theorem. Therefore
$$
\bigg\|\int_{\mathbb R^{\mathbb N}}(g(\bsy)-g_s(\bsy))\,\bsmu_\beta({\rm d}\bsy)\bigg\|_X\leq  \frac{C_{\rm init}}{(s^*)^{-\frac{2}{p}+1}}s^{-\frac{2}{p}+1}\quad\text{for all}~1\leq s\leq s^*.
$$
Thus it is enough to prove the claim for sufficiently large $s$. In what follows, we assume that $s> s^*$, let $G\in X'$ be arbitrary, and define
$$
F_G(\bsy):=\langle G,g(\bsy)\rangle_{X',X}\quad\text{for all}~\bsy\in U_{\boldsymbol\alpha}.
$$
Let $k$ be specified as in~\ref{assump2} (note that it always holds that $k\geq 2$). Then
$$
\partial^{\bsnu}F_G(\bsy)=\langle G,\partial^{\bsnu}g(\bsy)\rangle_{X',X}\quad\text{for all}~\bsnu\in\mathcal F_k~\text{and}~\bsy\in U_{\boldsymbol\alpha}
$$
and it follows from our assumptions that
$$
|\partial^{\bsnu}F_G(\bsy)|\leq \|G\|_{X'}\Theta_{|\bsnu|}\bsb^{\bsnu}\prod_{j\geq 1}{\rm e}^{\alpha_j|y_j|}\quad\text{for all}~\bsnu\in\mathcal F_k~\text{and}~\bsy\in U_{\boldsymbol\alpha}.
$$
Let $\bsy\in U_{\bsalpha}$ be arbitrary. Since $F_G$ is continuously differentiable up to order $k+1$, we can develop the Taylor expansion around the point $(\bsy_{\leq s},\mathbf 0)$, which yields
\begin{align*}
F_G(\bsy)&=F_G(\bsy_{\leq s},\mathbf 0)+\sum_{\ell=1}^k \sum_{\substack{|\bsnu|=\ell\\ \nu_j=0~\forall j\leq s}}\frac{\bsy^{\bsnu}}{\bsnu!}\partial^{\bsnu}F_G(\bsy_{\leq s},\mathbf 0)\\
&\quad +\sum_{\substack{|\bsnu|=k+1\\ \nu_j=0~\forall j\leq s}}\frac{k+1}{\bsnu!}\bsy^{\bsnu}\int_0^1 (1-\tau)^k\partial^{\bsnu}F_G(\bsy_{\leq s},\tau \bsy_{>s})\,{\rm d}\tau.
\end{align*}
Rearranging this equation and integrating both sides against the $\beta$-Gaussian product measure yields
\begin{align*}
&\int_{\mathbb R^{\mathbb N}}(F_G(\bsy)-F_G(\bsy_{\leq s},\mathbf 0))\,\bsmu_{\beta}({\mathrm d}\bsy)\\
&=\sum_{\ell=1}^k \sum_{\substack{|\bsnu|=\ell\\ \nu_j=0~\forall j\leq s}}\frac{1}{\bsnu!}\int_{\mathbb R^{\mathbb N}}\bsy^{\bsnu}\partial^{\bsnu}F_G(\bsy_{\leq s},\mathbf 0)\,\bsmu_{\beta}({\rm d}\bsy)\\
&\quad +\sum_{\substack{|\bsnu|=k+1\\ \nu_j=0~\forall j\leq s}}\frac{k+1}{\bsnu!}\int_{\mathbb R^{\mathbb N}}\int_0^1 (1-\tau)^k\bsy^{\bsnu}\partial^{\bsnu}F_G(\bsy_{\leq s},\tau \bsy_{>s})\,{\rm d}\tau\,\bsmu_\beta({\mathrm d}\bsy).
\end{align*}
If there exists a single component $\nu_k=1$ with $k>s$, then the summand in the first term vanishes since, by Lemma~\ref{lemma:infinite} and Fubini's theorem,
\begin{align*}
&\int_{\mathbb R^{\mathbb N}} \bsy^{\bsnu}\partial^{\bsnu}F_G(\bsy_{\leq s},\mathbf 0)\bsmu_{\beta}({\mathrm d}\bsy)\\
&=\!\bigg(\!\int_{\mathbb R^s}\!\partial^{\bsnu}F_G(\bsy_{\leq s},\mathbf 0)\!\prod_{j=1}^s\varphi_\beta(y_j)\,{\rm d}\bsy_{\leq s}\!\bigg)\!\underset{=0}{\underbrace{\bigg(\!\int_{\mathbb R}y_k\varphi_\beta(y_k)\,{\rm d}y_k\!\bigg)}}\!\bigg(\!\int_{\mathbb R^{\mathbb N}}\bsy^{\bsnu}\bsmu_{\beta}({\rm d}\bsy_{\{s+1:\infty\}\setminus\{k\}})\!\bigg)\\
&=0.
\end{align*}
Hence
\begin{align}
&\bigg|\int_{\mathbb R^{\mathbb N}}(F_G(\bsy)-F_G(\bsy_{\leq s},\mathbf 0))\,\bsmu_{\beta}({\mathrm d}\bsy)\bigg|\notag\\
&\leq \sum_{\ell=2}^k \sum_{\substack{|\bsnu|=\ell\\ \nu_j=0~\forall j\leq s\\ \nu_j\neq 1~\forall j>s}}\frac{1}{\bsnu!}\int_{\mathbb R^{\mathbb N}}|\bsy^{\bsnu}|\cdot|\partial^{\bsnu}F_G(\bsy_{\leq s},\mathbf 0)|\,\bsmu_{\beta}({\rm d}\bsy)\label{eq:term1}\\
&\quad +\sum_{\substack{|\bsnu|=k+1\\ \nu_j=0~\forall j\leq s}}\frac{k+1}{\bsnu!}\int_{\mathbb R^{\mathbb N}}\int_0^1 (1-\tau)^k|\bsy^{\bsnu}|\cdot|\partial^{\bsnu}F_G(\bsy_{\leq s},\tau \bsy_{>s})|\,{\rm d}\tau\,\bsmu_\beta({\mathrm d}\bsy).\label{eq:term2}
\end{align}

We start our estimation by splitting the terms in \eqref{eq:term1}:
\begin{align*}
&\int_{\mathbb R^{\mathbb N}}|\bsy^{\bsnu}|\cdot|\partial^{\bsnu}F_G(\bsy_{\leq s},\mathbf 0)|\,\bsmu_{\beta}({\rm d}\bsy) \leq \|G\|_{X'}\Theta_{|\bsnu|}\bsb^{\bsnu} \int_{\mathbb R^{\mathbb N}}|\bsy^{\bsnu}| \prod_{j=1}^{s} {\rm e}^{\alpha_j |y_j|}\bsmu_{\beta}(\mathrm d\bsy)\\ &\leq \|G\|_{X'}\Theta_{|\bsnu|}\bsb^{\bsnu} \int_{\mathbb R^{\mathbb N}}|\bsy^{\bsnu}| \prod_{j\geq 1} {\rm e}^{\alpha_j |y_j|} \bsmu_{\beta}(\mathrm d\bsy)\\
&=\|G\|_{X'}\Theta_{|\bsnu|}\bsb^{\bsnu} \bigg(\underbrace{\prod_{j\in {\rm supp}{(\bsnu)}} \int_{\mathbb{R}} |y_j|^{\nu_j} {\rm e}^{\alpha_j|y_j|}\varphi_\beta(y_j){\mathrm d}y_j\!}_{=:\rm term_1}\bigg)\!\bigg(\underbrace{\prod_{j\notin {\rm supp}{(\bsnu)}}  \int_{\mathbb{R}} {\rm e}^{\alpha_j|y_j|}\varphi_\beta(y_j){\mathrm d}y_j\!}_{=:\rm term_2}\bigg), 
\end{align*}
where the final step follows from Lemma~\ref{lemma:infinite} and Fubini's theorem. In order to bound ${\rm term_1}$, note that
\begin{align}\label{eq:Cdef}
C_{\alpha_j,\beta,\nu_j}:=\int_{\mathbb R}|y_j|^{\nu_j}{\rm e}^{\alpha_j|y_j|}\varphi_{\beta}(y_j)\,{\mathrm d}y_j<\infty
\end{align}
since we assumed that $\beta\geq1$ and $\alpha_j<1$ in the case $\beta =1$. We define an auxiliary constant $A_{\alpha_j,\beta,\ell}:=\max_{1\leq k\leq \ell}C_{\alpha_j,\beta,k}$. Clearly, $C_{\alpha_j,\beta,\nu_j}\leq A_{\alpha_j,\beta,|\bsnu|}<\infty$ and $C_{\alpha,\beta,\nu}\leq C_{\alpha',\beta,\nu}$ whenever $\alpha\leq \alpha'$ with $\beta,\nu$ fixed. In particular,
\begin{align*}
{\rm term_1} &= \prod_{j\in {\rm supp}{(\bsnu)}} C_{\alpha_j,\beta,\nu_j} \leq  \prod_{j\in {\rm supp}{(\bsnu)}} C_{\|\bsalpha\|_{\infty},\beta,\nu_j} \leq \prod_{j\in {\rm supp}{(\bsnu)}} A_{\|\bsalpha\|_{\infty},\beta,|\bsnu|}\\ &\leq \max{\{1,A_{\|\bsalpha\|_{\infty},\beta,|\bsnu|}\}}^{|\bsnu|},
\end{align*}
where $\|\bsalpha\|_\infty:=\sup_{j\geq 1}|\alpha_j|$ is finite since $\bsalpha\in\ell^1(\mathbb N)$ by assumption. To bound ${\rm term_2}$, we note that there is an index $j'\in\mathbb N$ such that $\alpha_j\leq\frac12$ for all $j>j'$. Hence
\begin{align*}
{\rm term}_2&=\bigg(\prod_{\substack{j\not\in {\rm supp}{(\bsnu)}\\ 1\leq j\leq j'}}  \int_{\mathbb{R}} {\rm e}^{\alpha_j|y_j|}\varphi_\beta(y_j)\,{\mathrm d}y_j\bigg)\bigg(\prod_{\substack{j\not\in {\rm supp}{(\bsnu)}\\ j>j'}} \int_{\mathbb{R}} {\rm e}^{\alpha_j|y_j|}\varphi_\beta(y_j)\,{\mathrm d}y_j\bigg)\\
&\leq \max\{1,C_{\|\bsalpha\|_\infty,\beta,0}\}^{j'}\bigg(\prod_{\substack{j\not\in {\rm supp}{(\bsnu)}\\ j>j'}} \int_{\mathbb{R}} {\rm e}^{\alpha_j|y_j|}\varphi_\beta(y_j)\,{\mathrm d}y_j\bigg),
\end{align*}
using a similar argument as before. In order to ensure that the remaining factor is finite, we argue similarly to~\cite[Proposition~3.2]{hks21}.

Let us first consider the case $\beta >1$. Young's inequality states, for all $x,y\geq 0$ and $\theta\in(0,1)$, that
$$
xy=x^{\theta}x^{1-\theta}y\leq \frac{\beta-1}{\beta}x^{\theta\frac{\beta}{\beta-1}}+\frac{1}{\beta}x^{(1-\theta)\beta}y^\beta,
$$
where $\beta>1$. The special choice $\theta=\frac{\beta-1}{\beta}$ yields
$$
xy\leq \frac{\beta-1}{\beta}x+\frac{1}{\beta}xy^\beta\quad\text{for all}~x,y\geq 0.
$$
Thereby
\begin{align*}
\int_{\mathbb R}{\rm e}^{\alpha_j|y_j|}\varphi_\beta(y_j)\,{\rm d}y_j&\leq {\rm e}^{\frac{\beta-1}{\beta}\alpha_j}\int_{\mathbb R}{\rm e}^{\frac{\alpha_j}{\beta}|y_j|^\beta}\varphi_\beta(y_j)\,{\rm d}y_j\\
&=\frac{{\rm e}^{\frac{\beta-1}{\beta}\alpha_j}}{2\beta^{\frac{1}{\beta}}\Gamma(1+\frac{1}{\beta})}\int_{\mathbb R}{\rm e}^{-(1-\alpha_j)\frac{|y_j|^\beta}{\beta}}\,{\rm d}y_j
=\frac{{\rm e}^{\frac{\beta-1}{\beta}\alpha_j}}{(1-\alpha_j)^{\frac{1}{\beta}}},
\end{align*}
where we used $\int_{\mathbb R}{\rm e}^{-(1-\alpha_j)\frac{|y_j|^\beta}{\beta}}\,{\rm d}y_j=2\beta^{\frac{1}{\beta}}(1-\alpha_j)^{-\frac{1}{\beta}}\Gamma(1+\frac{1}{\beta})$. Furthermore, since $$\frac{1}{1-x}=1+\frac{x}{1-x}\leq \exp\bigg(\frac{x}{1-x}\bigg)\quad\text{for all}~x\in[0,1),$$
we obtain
\begin{align*}
\int_{\mathbb R}{\rm e}^{\alpha_j|y_j|}\varphi_\beta(y_j)\,{\rm d}y_j\leq\exp\bigg(\frac{\beta\!-\!1}{\beta}\alpha_j\bigg)\exp\bigg(\frac{1}{\beta}\frac{\alpha_j}{1\!-\!\alpha_j}\bigg)\leq \exp\bigg(\frac{\beta\!+\!1}{\beta}\alpha_j\bigg)\leq \exp(2\alpha_j)
\end{align*}
since we assumed $\alpha_j\leq\frac12$ for all $j>j'$. Therefore
$$
\prod_{\substack{j\notin {\rm supp}(\bsnu)\\ j>j'}}\int_{\mathbb R}{\rm e}^{\alpha_j|y_j|}\varphi_\beta(y_j)\,{\rm d}y_j\leq\exp\bigg(2\!\!\sum_{\substack{j\notin {\rm supp}(\bsnu)\\ j>j'}}\!\alpha_j\bigg) \leq\exp\bigg(2\sum_{j\geq 1}\alpha_j\bigg)=:\widetilde{C},
$$
where $\widetilde{C}<\infty$ since $\bsalpha\in\ell^1(\mathbb N)$.
The special case $\beta=1$ follows since we assumed $\alpha_j<1$ for all $j\geq 1$ and it holds that
$$
\int_{\mathbb R}{\rm e}^{\alpha_j|y_j|}\varphi_1(y_j)\,{\rm d}y_j=\frac{1}{1-\alpha_j}\leq \exp\bigg(\frac{\alpha_j}{1-\alpha_j}\bigg)\leq \exp(2\alpha_j)
$$
for $j>j'$. Thus
$$
\prod_{\substack{j\notin {\rm supp}(\bsnu)\\ j>j'}}\int_{\mathbb R}{\rm e}^{\alpha_j|y_j|}\varphi_1(y_j)\,{\rm d}y_j\leq \exp\bigg(2\sum_{\substack{j\notin {\rm supp}(\bsnu)\\ j>j'}}\alpha_j\bigg)\leq \exp\bigg(2\sum_{j\geq 1}\alpha_j\bigg)=\widetilde{C}<\infty
$$
since $\bsalpha\in\ell^1(\mathbb N)$.
Combining the estimates for ${\rm term_1}$ and ${\rm term_2}$ gives
\begin{align*}
    &\int_{\mathbb R^{\mathbb N}}|\bsy^{\bsnu}|\cdot|\partial^{\bsnu}F_G(\bsy_{\leq s},\mathbf 0)|\,\bsmu_{\beta}({\rm d}\bsy)\\ &\leq \widetilde{C}\|G\|_{X'}\max\{1,C_{\|\bsalpha\|_\infty,\beta,0}\}^{j'}\Theta_{|\bsnu|}\bsb^{\bsnu}  \max{\{1,A_{\|\alpha\|_{\infty},\beta,|\bsnu|}\}}^{|\bsnu|},
\end{align*}
where $C_{\|\bsalpha\|_\infty,\beta,0}$ is defined in \eqref{eq:Cdef}.
Similarly we split the terms in \eqref{eq:term2},
\begin{align*}
    &\int_{\mathbb R^{\mathbb N}}\int_0^1 (1-\tau)^k|\bsy^{\bsnu}|\cdot|\partial^{\bsnu}F_G(\bsy_{\leq s},\tau \bsy_{>s})|\,{\rm d}\tau\,\bsmu_\beta({\mathrm d}\bsy)\\
    &\leq \|G\|_{X'}\Theta_{|\bsnu|}\bsb^{\bsnu} \bigg( \prod_{j \in {\rm supp}(\bsnu)} \int_{\mathbb R} |y_j|^{\nu_j} {\rm e}^{\alpha_j |y_j|} \varphi_{\beta}(y_j) \mathrm dy_j\bigg)
    \bigg(\prod_{j \notin {\rm supp}(\bsnu)} \int_{\mathbb R} {\rm e}^{\alpha_j |y_j|} \varphi_{\beta}(y_j) \mathrm dy_j\bigg)\\
    &\leq \widetilde{C}\|G\|_{X'}\Theta_{|\bsnu|}\max\{1,C_{\|\bsalpha\|_\infty,\beta,0}\}^{j'}\bsb^{\bsnu} \max{\{1,A_{\|\bsalpha\|_{\infty},\beta,|\bsnu|}\}}^{|\bsnu|},
\end{align*}
where we used 
\begin{align*}
    |\partial^{\bsnu}F_G(\bsy_{\leq s},\tau \bsy_{>s})|&\leq \|G\|_{X'}\Theta_{|\bsnu|}\bsb^{\bsnu} \bigg(\prod_{j=1}^s {\rm e}^{\alpha_j|y_j|}\bigg)\bigg(\prod_{j>s}{\rm e}^{\tau\alpha_j|y_j|}\bigg)\\ &\leq \|G\|_{X'}\Theta_{|\bsnu|}\bsb^{\bsnu} \bigg(\prod_{j\geq1} {\rm e}^{\alpha_j|y_j|}\bigg).
\end{align*}
These inequalities allow us to estimate
\begin{align}
&\bigg|\int_{\mathbb R^{\mathbb N}}(F_G(\bsy)-F_G(\bsy_{\leq s},\mathbf 0))\,\bsmu_{\beta}({\mathrm d}\bsy)\bigg|\notag\\
&\leq \widetilde{C}\|G\|_{X'}\max\{1,C_{\|\bsalpha\|_\infty,\beta,0}\}^{j'}\sum_{\ell=2}^k \sum_{\substack{|\bsnu|=\ell\\ \nu_j=0~\forall j\leq s\\ \nu_j\neq 1~\forall j>s}}\frac{\Theta_{\ell}}{\bsnu!}\bsb^{\bsnu}  \max\{1,A_{\|\bsalpha\|_\infty,\beta,\ell}\}^{\ell}\notag\\
&\quad +\widetilde{C}\|G\|_{X'}\max\{1,C_{\|\bsalpha\|_\infty,\beta,0}\}^{j'}\sum_{\substack{|\bsnu|=k+1\\ \nu_j=0~\forall j\leq s}}\frac{k+1}{\bsnu!}\Theta_{k+1}\bsb^{\bsnu}  \max\{1,A_{\|\bsalpha\|_\infty,\beta,k+1}\}^{k+1}\notag\\
&\leq \widetilde{C}\|G\|_{X'}\max\{1,C_{\|\bsalpha\|_\infty,\beta,0}\}^{j'} \big(\max_{2\leq\ell\leq k}(\Theta_{\ell}\max\{1,A_{\|\bsalpha\|_\infty,\beta,\ell}\}^\ell)\big)\sum_{\ell=2}^k \sum_{\substack{|\bsnu|=\ell\\ \nu_j=0~\forall j\leq s\\ \nu_j\neq 1~\forall j>s}}\bsb^{\bsnu}\label{eq:gantner}\\
&\quad +\widetilde{C}\|G\|_{X'}\max\{1,C_{\|\bsalpha\|_\infty,\beta,0}\}^{j'}\Theta_{k+1}(k+1)\max\{1,A_{\|\bsalpha\|_\infty,\beta,k+1}\}^{k+1} \sum_{\substack{|\bsnu|=k+1\\ \nu_j=0~\forall j\leq s}}\bsb^{\bsnu}.\label{eq:kss}
\end{align}

To bound the term~\eqref{eq:gantner}, we argue similarly to \cite[Theorem~1]{gantner} by noting that it follows from our definition of $s^*$ that
$$
\bigg(\sum_{k>s}b_k\bigg)^2\leq \frac14\quad\text{and} \quad b_j\leq\frac{1}{2}\quad\text{for all}~j>s.
$$
This leads us to estimate
\begin{align*}
&\sum_{\ell=2}^k \sum_{\substack{|\bsnu|=\ell\\ \nu_j=0~\forall j\leq s\\ \nu_j\neq 1~\forall j>s}}\bsb^{\bsnu}=\sum_{\substack{2\leq |\bsnu|\leq k\\ \nu_j=0~\forall j\leq s\\ \nu_j\neq 1~\forall j>s}}\bsb^{\bsnu}\leq \sum_{\substack{0\neq |\bsnu|_\infty\leq k\\ \nu_j=0~\forall j\leq s\\ \nu_j\neq 1~\forall j>s}}\bsb^{\bsnu}=-1+\prod_{j>s}\bigg(1+\sum_{\ell=2}^k b_j^{\ell}\bigg)\\
&=-1+\prod_{j>s}\bigg(1+\frac{1-b_j^{k-1}}{1-b_j}b_j^2\bigg)\leq -1+\prod_{j>s}\bigg(1+2b_j^2\bigg)\leq -1+\exp\bigg(2\sum_{j>s}b_j^2\bigg)\\
&\leq 2({\rm e}-1)\sum_{j>s}b_j^2,
\end{align*}
where the final inequality is a consequence of Bernoulli's inequality $(1+x)^r\leq 1+rx$ for all $0\leq r\leq 1$ and $x\geq-1$. It is an immediate consequence of Lemma~\ref{lemma:stechkin} that
\begin{align}
\sum_{j>s}b_j^2\leq s^{-\frac{2}{p}+1}\bigg(\sum_{j>s}b_j^p\bigg)^{\frac{2}{p}}\label{eq:stechkin1}
\end{align}
since $\bsb$ was assumed to be a nonincreasing sequence such that $\bsb\in\ell^p(\mathbb N)$ for some $p\in(0,1)$.

We estimate the term~\eqref{eq:kss} similarly to the approach taken in~\cite[Theorem~4.1]{GGKSS2019}. By the trivial bound $\frac{|\bsnu|!}{\bsnu!}\geq 1$ and the multinomial theorem, we obtain
\begin{align}
\sum_{\substack{|\bsnu|=k+1\\ \nu_j=0~\forall j\leq s}}\bsb^{\bsnu}\leq \sum_{\substack{|\bsnu|=k+1\\ \nu_j=0~\forall j\leq s}}\frac{|\bsnu|!}{\bsnu!}\bsb^{\bsnu}=\bigg(\sum_{j>s}b_j\bigg)^{k+1}\leq s^{(-\frac{1}{p}+1)(k+1)}\bigg(\sum_{j>s}b_j^p\bigg)^{(k+1)/p},\label{eq:stechkin2}
\end{align}
where the final inequality follows immediately from Lemma~\ref{lemma:stechkin} and our assumption that $\bsb$ is a nonincreasing sequence such that $\bsb\in \ell^p(\mathbb N)$ for some $p\in(0,1)$.

Putting the inequalities~\eqref{eq:stechkin1} and~\eqref{eq:stechkin2} together and utilizing $k=\lceil \frac{1}{1-p}\rceil$ (see \ref{assump2}) we obtain
$$
\bigg|\int_{\mathbb R^{\mathbb N}}(F_G(\bsy)-F_G(\bsy_{\leq s},\mathbf 0))\,\bsmu_{\beta}({\mathrm d}\bsy)\bigg|\leq C\|G\|_{X'}(s^{-\frac{2}{p}+1}+s^{(-\frac{1}{p}+1)(k+1)})\!\leq\! C\|G\|_{X'}s^{-\frac{2}{p}+1}
$$
for some constant $C>0$ independent of $s$. Finally, by recalling that $F_G(\bsy)=\langle G,g(\bsy)\rangle_{X',X}$ and $G\in X'$ was arbitrary, we can take the supremum over $\{G\in X':\|G\|_{X'}\leq 1\}$ to obtain
\begin{align*}
    &\sup_{G\in X':\|G\|_{X'}\leq1}\bigg|\int_{\mathbb R^{\mathbb N}}(F_G(\bsy)-F_G(\bsy_{\leq s},\mathbf 0))\,\bsmu_{\beta}({\mathrm d}\bsy)\bigg|\\
    &=\sup_{G\in X':\|G\|_{X'}\leq1}\bigg|\int_{\mathbb R^{\mathbb N}}\langle G,g(\bsy)-g(\bsy_{\leq s},\mathbf 0)\rangle_{X',X}\,\bsmu_{\beta}({\mathrm d}\bsy)\bigg|\\
    &=\sup_{G\in X':\|G\|_{X'}\leq1}\bigg|\bigg\langle G,\int_{\mathbb R^{\mathbb N}}(g(\bsy)-g(\bsy_{\leq s},\mathbf 0))\,\bsmu_{\beta}({\mathrm d}\bsy)\bigg\rangle_{X',X}\bigg|\\
    &=\bigg\|\int_{\mathbb R^{\mathbb N}}(g(\bsy)-g_s(\bsy))\,\bsmu_\beta({\mathrm d}\bsy)\bigg\|_X \leq Cs^{-\frac{2}{p}+1}
\end{align*}
as desired.\quad\endproof

In general the constant~$C$ in the statement of~Theorem~\ref{thm:main} depends on~$\beta \ge 1$. Under additional assumptions on the sequence~$\bsalpha$ it is possible to bound the constant independently of~$\beta$.
\begin{corollary}
    Let~$\|\bsalpha\|_{\infty} <1$. Then the constant~$C$ in~Theorem~\ref{thm:main} can be bounded independently of~$\beta\geq 1$.
\end{corollary}
\proof
We note that~$\beta$ enters the constant in~Theorem~\ref{thm:main} only through~$C_{\|\bsalpha\|_\infty,\beta,\nu_j}$ as defined in~\eqref{eq:Cdef}. Thus, it is sufficient to find an upper bound for
\begin{align*}
    C_{\alpha,\beta,\nu} = \int_\bbR |y|^{\nu} {\rm e}^{\alpha|y|}\varphi_\beta(y) \,\mathrm dy
\end{align*}
with~$\alpha = \|\bsalpha\|_\infty$ and $\nu = \nu_j$ as well as $\nu = 0$, independently of~$\beta$.
By~\cite[Corollary 1]{Dytso18},
$$
\int_\bbR |y|^{\nu} \varphi_\beta(y) \,\mathrm dy = \frac{\beta^{\frac{\nu}{\beta}}}{\Gamma(\frac{1}{\beta})}\Gamma\bigg(\frac{\nu+1}{\beta}\bigg) \le \Gamma(\nu+1)
= \int_\bbR |y|^{\nu} \varphi_1(y) \,\mathrm dy
$$
for all $\beta \ge 1$ and $\nu\geq 0$. We conclude that
$$
\int_\bbR |y|^{\nu} {\rm e}^{\alpha|y|}\varphi_\beta(y) \,\mathrm dy  \le \int_\bbR |y|^{\nu} {\rm e}^{\alpha|y|}\varphi_1(y) \,\mathrm dy
$$
holds with $\alpha = 0$. Next, we show that this inequality is valid for $\alpha\ge 0$. This follows by observing that
\begin{align*}
\frac{\partial}{\partial\alpha}\int_{\mathbb R}|y|^{\nu}{\rm e}^{\alpha|y|}\varphi_\beta(y)\,{\rm d}y
&=\int_{\mathbb R}|y|^{\nu+1}{\rm e}^{\alpha|y|}\varphi_\beta(y)\,{\rm d}y\\
&\leq\int_{\mathbb R}|y|^{\nu+1}{\rm e}^{\alpha|y|}\varphi_1(y)\,{\rm d}y\\
&=\frac{\partial}{\partial\alpha}\int_{\mathbb R}|y|^{\nu}{\rm e}^{\alpha|y|}\varphi_1(y)\,{\rm d}y.
\end{align*}
Altogether we obtain that
\begin{align*}
    C_{\alpha,\beta,\nu} \le C_{\alpha,1,\nu} < \infty
\end{align*}
since~$\alpha<1$.\quad\endproof

We also state the corresponding dimension truncation result in the uniform case formally corresponding to $\beta=\infty$.
\begin{theorem}\label{thm:main2}
Let $g(\bsy)\in X$ for all $\bsy\in [-1,1]^{\mathbb N}$. Suppose that assumptions {\rm\ref{assump1pr}} and {\rm\ref{assump2pr}} hold. Then
$$
\bigg\|\int_{\mathbb [-1,1]^{\mathbb N}}(g(\bsy)-g_s(\bsy))\,\bsgamma({\mathrm d}\bsy)\bigg\|_X\leq Cs^{-\frac{2}{p}+1},
$$
where the constant $C>0$ is independent of the dimension $s$.

Let $G\in X'$ be arbitrary. Then
$$
\bigg|\int_{\mathbb [-1,1]^{\mathbb N}}G(g(\bsy)-g_s(\bsy))\,\bsgamma({\mathrm d}\bsy)\bigg|\leq C\|G\|_{X'}s^{-\frac{2}{p}+1},
$$
where the constant $C>0$ is as above.
\end{theorem}
\proof The steps are completely analogous to Theorem~\ref{thm:main} in the special case $\alpha_j=0$ for all $j\geq 1$ and by restricting the domain of integration to $[-1,1]^{\mathbb N}$. In the special case $\Theta_{|\bsnu|}=(|\bsnu|+r_1)!$, $r_1\in\mathbb N_0$, the proof works as in~\cite[Theorem~6.2]{guth22}.\quad\endproof

\emph{Remark.} The conditions~\ref{assump2} and~\ref{assump2pr} are formulated as sufficient conditions. The form in which the regularity bounds are postulated in~\ref{assump2} and~\ref{assump2pr} is an important ingredient for the Taylor series argument. However, it is known that~\ref{assump2pr} is not a necessary condition: an example is given in~\cite[Lemma 2.4]{kks20}, where the authors obtain the dimension truncation rate $\mathcal O(s^{-\frac{2}{p}+1})$ for a problem which satisfies a more general parametric regularity bound than~\ref{assump2pr}.

\section{Application to parametric elliptic PDEs}\label{sec:pdeappl}
In this section, we illustrate how to apply the main dimension truncation results proved in Section~\ref{sec:dimtrunc} to parametric elliptic PDE model problems. We consider uniform and affine as well as lognormal parameterizations of the input random field. The rate we obtain for the uniform and affine model coincides with the well-known rate in the literature~\cite{gantner} and is not a new result, however, we present it for completeness. Remarkably, the dimension truncation rate we obtain for the lognormal model using our method improves the rates in the existing literature (cf., e.g.,~\cite{log,kuonuyenssurvey}). Finally, we give an example on how our results can be applied to assess dimension truncation rates corresponding to PDE solutions composed with nonlinear quantities of interest.

We consider the problem of finding $u:D\times U\to \bbR$ such that
\begin{align}
\begin{split}
-\nabla\cdot(a(\bsx,\bsy)\nabla u(\bsx,\bsy))&=f(\bsx), \quad\bsx\in D,~\bsy\in U,\\
u(\bsx,\bsy)&=0, \quad\quad\bsx\in\partial D,~\bsy\in U,
\end{split}\label{eq:pdeproblem}
\end{align}
in some bounded Lipschitz domain $D\subset\bbR^d$, $d\in \{1,2,3\}$, for some given source term $f\!:D\to\mathbb R$ and diffusion coefficient $a\!:D\times U\to \mathbb R$. The parameter set $U$ is assumed to be a nonempty subset of $\mathbb R^{\mathbb N}$.%

The relevant function spaces for the elliptic PDE problem~\eqref{eq:pdeproblem} are $X:=H_0^1(D)$ and its dual $X'=H^{-1}(D)$, which we understand with respect to the pivot space $H:=L^2(D)$. The space $H$ is identified with its own dual and we set
$$
\|v\|_X:=\|\nabla v\|_{H}\quad\text{for}~v\in X.
$$
The weak formulation of~\eqref{eq:pdeproblem} is to find, for all $\bsy\in U$, a solution $u(\cdot,\bsy)\in X$ such that
\begin{align}
\int_Da(\bsx,\bsy)\nabla u(\bsx,\bsy)\cdot\nabla v(\bsx)\,{\rm d}\bsx=\langle f,v\rangle_{X',X}\quad\text{for all}~v\in X,\label{eq:pdeproblem2}
\end{align}
where $f\in X'$ and $\langle\cdot,\cdot\rangle_{X',X}$ denotes the duality pairing between elements of $X'$ and $X$.

The following lemma collects basic, well-known results about the existence of a unique solution to~\eqref{eq:pdeproblem2} (Lax--Milgram lemma), the continuity of the PDE solution with respect to the right-hand side of~\eqref{eq:pdeproblem2} (\emph{a priori} bound), and the continuity of the PDE solution with respect to the diffusion coefficient (the second Strang lemma).

\begin{lemma}\label{lemma:strang} Let $D\subset \mathbb R^d$, $d\in\{1,2,3\}$, be a bounded Lipschitz domain, $\varnothing\neq U\subseteq\mathbb R^{\mathbb N}$, $f\in X'$, and suppose that there exist $a_{\min}(\bsy):=\min_{\bsx\in\overline{D}}a(\bsx,\bsy)\in L^\infty_+(D)$ and $a_{\max}(\bsy):=\max_{\bsx\in\overline{D}}a(\bsx,\bsy)\in L^\infty_+(D)$ such that 
\begin{align}
0<a_{\min}(\bsy)\leq a(\bsx,\bsy)\leq a_{\max}(\bsy)<\infty\quad\text{for all}~\bsx\in D~\text{and}~\bsy\in U,\label{eq:unifelliptic}
\end{align}
where $a(\cdot,\bsy)\in L^\infty_+(D)$, $\bsy\in U$, is the diffusion coefficient in~\eqref{eq:pdeproblem2}. We define $a_s(\cdot,\bsy):=a(\cdot,(\bsy_{\leq s},\mathbf 0))$, $a_{\min}^s(\bsy):=a_{\min}(\bsy_{\leq s},\mathbf 0)$, and $u_s(\cdot,\bsy):=u(\cdot,(\bsy_{\leq s},\mathbf 0))$. Then there exists a unique solution to~\eqref{eq:pdeproblem2} such that
\begin{align}
\|u(\cdot,\bsy)\|_X\leq \frac{\|f\|_{X'}}{a_{\min}(\bsy)}\quad\text{for all}~\bsy\in U\label{eq:apriori}
\end{align}
and
\begin{align}\label{eq:strang}
\|u(\cdot,\bsy)-u_s(\cdot,\bsy)\|_X\leq \frac{1}{a_{\min}(\bsy)a_{\min}^s(\bsy)}\|a(\cdot,\bsy)-a_s(\cdot,\bsy)\|_{L^\infty(D)}\|f\|_{X'}\quad\!\text{for all}~\bsy\in U.
\end{align}
\end{lemma}\proof For proofs of these results we refer to~\cite{matern1}.\quad\endproof

Studies in uncertainty quantification for PDEs typically consider one of the following two models for the input random field.
\begin{enumerate}[label=(M\arabic*)]
    \item The diffusion coefficient is parameterized by%
    $$
    a(\bsx,\bsy)=a_0(\bsx)\exp\bigg(\sum_{j=1}^\infty y_j\psi_j(\bsx)\bigg),\quad y_j\in \mathbb R,%
    $$
    for $a_0 \in L^{\infty}_+(D)$, $\psi_j\in L^\infty(D)$ for all $j\geq 1$ with $(\|\psi_j\|_{L^\infty})\in\ell^p(\mathbb N)$ for some $p\in(0,1)$, and $U = \bbR^\bbN$. Hence, it follows that
    \begin{align*}
        0<a_{\min}(\bsy) \leq a(\bsx,\bsy) \leq a_{\max}(\bsy) < \infty\quad \text{for all }\bsx \in D \text{ and }\bsy \in U,
    \end{align*}
    where $a_{\min}(\bsy) = \min_{x\in \overline{D}} a(\bsx,\bsy)$ and $a_{\max}(\bsy) = \max_{x\in \overline{D}} a(\bsx,\bsy)$.\label{model1}
        \item The diffusion coefficient is parameterized by
    $$
    a(\bsx,\bsy)=a_0(\bsx)+\sum_{j=1}^\infty y_j\psi_j(\bsx),\quad y_j\in [-1,1],
    $$  
    for $a_0\in L^{\infty}(D)$, $\psi_j \in L^{\infty}(D)$ for all $j\geq 1$ with $(\|\psi_j\|_{L^\infty})\in\ell^p(\mathbb N)$ for some $p\in(0,1)$, and $U = [-1,1]^\bbN$, such that
    \begin{align*}
        0<a_{\min} \leq a(\bsx,\bsy) \leq a_{\max} < \infty\quad \text{for all }\bsx \in D \text{ and }\bsy \in U,
    \end{align*}
    for some constants $a_{\max}\geq a_{\min}>0$ independent of $\bsx\in D$ and $\bsy\in U$.\label{model2}
\end{enumerate}
Let $\bsb:=(b_j)_{j\geq 1}$ with $b_j:=\|\psi_j\|_{L^\infty}$. In addition, we assume that $b_1\geq b_2\geq\cdots$.

Recall that $u_s(\cdot,\bsy)=u(\cdot,(\bsy_{\leq s},\mathbf 0))$ for $\bsy\in U$. In the context of high-dimensional numerical integration, it is germane in the setting~\ref{model1} to quantify the dimension truncation error
$$
\bigg\|\int_{U}(u(\cdot,\bsy)-u_s(\cdot,\bsy))\,\bsmu_{\beta}({\rm d}\bsy)\bigg\|_X,\quad\text{with}~U=\mathbb R^{\mathbb N},
$$
and in the setting~\ref{model2}, the dimension truncation error
$$
\bigg\|\int_{U}(u(\cdot,\bsy)-u_s(\cdot,\bsy))\,\bsgamma({\rm d}\bsy)\bigg\|_X,\quad\text{with}~U=[-1,1]^{\mathbb N}.
$$
\paragraph{Lognormal model and its generalizations.} The model~\ref{model1} is the \emph{lognormal model} (cf., e.g.,~\cite{gittelson10,log,log2,log3,log4,log5,schwabtodor}) when the uncertain parameter $\bsy\in \mathbb R^{\mathbb N}$ is endowed with the $\beta$-Gaussian probability measure with $\beta=2$ and $\alpha_j=b_j$ for all $j\geq 1$. However, the dimension truncation analysis in Section~\ref{sec:dimtrunc} covers the more general setting where we have arbitrary $\bsalpha\in \ell^1(\mathbb N)$ and either $\beta\in(1,\infty)$ or $\beta=1$ with $\alpha_j<1$ for all $j\geq 1$. We remark that the latter case corresponds to random variables distributed according to the Laplace distribution. By~\eqref{eq:strang}, condition~\ref{assump1} holds because
$$
\|a(\cdot,\bsy)-a_s(\cdot,\bsy)\|_{L^\infty(D)}\xrightarrow{s\to\infty}0\quad\text{for all}~\bsy\in U_{{\bsb}}.
$$
On the other hand, condition~\ref{assump2} holds due to the well-known (see, e.g., \cite[Theorem 14]{log}) parametric regularity bound%
$$
\|\partial^{\bsnu}u(\cdot,\bsy)\|_X\leq \frac{\|f\|_{X'}}{\min_{\bsx\in \overline{D}}{a_0(\bsx)}}\frac{|\bsnu|!}{(\log 2)^{|\bsnu|}}\bsb^\bsnu \prod_{j\geq 1}{\rm e}^{b_j|y_j|}\quad\text{for all}~\bsy\in U_{\bsb},~\bsnu\in \mathcal F.
$$
Especially, this corresponds to our setting with the special choice $\alpha_j=b_j$ for all $j\geq 1$. 
By applying Theorem~\ref{thm:main} to $u$, we obtain that
$$
\bigg\|\int_{\mathbb R^{\mathbb N}}(u(\cdot,\bsy)-u_s(\cdot,\bsy))\,\bsmu_{\beta}({\rm d}\bsy)\bigg\|_X=\mathcal O(s^{-\frac{2}{p}+1}),
$$
where the implied coefficient is independent of the dimension $s$.

Let $D$ be a bounded polyhedron and suppose that $\{X_h\}_h$ is a family of conforming finite element subspaces $X_h\subset X$, indexed by the mesh size $h>0$. Let $u_h(\cdot,\bsy)\in X_h$ and $u_{s,h}(\cdot,\bsy)\in X_h$ denote the finite element discretized solutions corresponding to $u(\cdot,\bsy)$ and $u_s(\cdot,\bsy)$, respectively. Then it also holds that
$$
\bigg\|\int_{\mathbb R^{\mathbb N}}(u_h(\cdot,\bsy)-u_{s,h}(\cdot,\bsy))\,\bsmu_{\beta}({\rm d}\bsy)\bigg\|_X=\mathcal O(s^{-\frac{2}{p}+1}),
$$
where the implied coefficient is again independent of the dimension $s$.
\paragraph{Uniform and affine model.} The model~\ref{model2} is known as the \emph{uniform and affine model} (cf., e.g.,~\cite{cds10,dicklegiaschwab,spodpaper14,ghs18,schwab_gittelson_2011,kss12,kssmultilevel,schwab13}) when the uncertain parameter $\bsy\in [-1,1]^{\mathbb N}$ is endowed with the uniform probability measure. By~\eqref{eq:strang}, condition~\ref{assump1pr} holds since
$$
\|a(\cdot,\bsy)-a_s(\cdot,\bsy)\|_{L^\infty(D)}\xrightarrow{s\to\infty}0\quad \text{for all}~\bsy\in[-1,1]^{\mathbb N}.
$$
Moreover, condition~\ref{assump2pr} holds due to the well-known (see, e.g., \cite{cds10}) parametric regularity bound
\begin{align}
\|\partial^{\bsnu}u(\cdot,\bsy)\|_X \leq \frac{\|f\|_{X'}}{a_{\min}^{|\bsnu|+1}}|\bsnu|!\bsb^{\bsnu}\quad\text{for all}~\bsy\in [-1,1]^{\mathbb N},~\bsnu\in\mathcal F.\label{eq:unifregu}
\end{align}
It follows from applying Theorem~\ref{thm:main2} to $u$ that
$$
\bigg\|\int_{[-1,1]^{\mathbb N}}(u(\cdot,\bsy)-u_s(\cdot,\bsy))\,\bsgamma({\rm d}\bsy)\bigg\|_X=\mathcal O(s^{-\frac{2}{p}+1}),
$$
where the implied coefficient is independent of the dimension $s$. The same result holds if $u$ and $u_s$ are replaced by finite element solutions belonging to a conforming finite element subspace of $X$.

Finally, we present an example illustrating how our results can be applied to nonlinear quantities of interest of the PDE response.

\begin{example}\label{ex:nonlin}\rm
Consider the uniform and affine model (M2) with $U=[-1,1]^{\mathbb N}$ and suppose that we are interested in analyzing
\begin{align}
\bigg|\int_{[-1,1]^{\mathbb N}}(G_{\rm nl}(u(\cdot,\bsy))-G_{\rm nl}(u_s(\cdot,\bsy)))\,\bsgamma({\mathrm d}\bsy)\bigg|,\label{eq:nonlinearex}
\end{align}
where $u$ and $u_s$ denote the parametric PDE solution and its dimension truncation in $X=H^1_0(D)$, respectively. Suppose that the quantity of interest is the nonlinear functional
\begin{align}
G_{\rm nl}(v):=\int_D v(\bsx)^2\,{\rm d}\bsx,\quad v\in X.\label{eq:nonlinearqoi}
\end{align}
Let $\bsnu\in\mathcal F$ and $\bsy\in[-1,1]^{\mathbb N}$. It follows by an application of the Leibniz product rule and the regularity bound~\eqref{eq:unifregu} that
\begin{align*}
\partial^{\bsnu}G_{\rm nl}(u(\cdot,\bsy))&=\sum_{\bsm\leq \bsnu}\binom{\bsnu}{\bsm}\int_D\partial^{\bsm}u(\bsx,\bsy)\cdot \partial^{\bsnu-\bsm}u(\bsx,\bsy)\,{\rm d}\bsx\\
&\leq \sum_{\bsm\leq \bsnu}\binom{\bsnu}{\bsm}\|\partial^{\bsm}u(\cdot,\bsy)\|_{L^2(D)}\|\partial^{\bsnu-\bsm}u(\cdot,\bsy)\|_{L^2(D)}\\
&\leq C_P^2\sum_{\bsm\leq \bsnu}\binom{\bsnu}{\bsm}\|\partial^{\bsm}u(\cdot,\bsy)\|_{X}\|\partial^{\bsnu-\bsm}u(\cdot,\bsy)\|_{X}\\
&\leq\frac{C_P^2\|f\|_{X'}^2}{a_{\min}^{|\bsnu|+2}}\bsb^{\bsnu}\sum_{\bsm\leq\bsnu}\binom{\bsnu}{\bsm}|\bsm|!\,|\bsnu-\bsm|!\\
&=\frac{C_P^2\|f\|_{X'}^2}{a_{\min}^{|\bsnu|+2}}\bsb^{\bsnu}\sum_{\ell=0}^{|\bsnu|}\ell!\,(|\bsnu|-\ell)!\sum_{\substack{\bsm\leq \bsnu\\|\bsm|=\ell}}\binom{\bsnu}{\bsm}\\
&=\frac{C_P^2\|f\|_{X'}^2}{a_{\min}^{|\bsnu|+2}}(|\bsnu|+1)!\,\bsb^{\bsnu},
\end{align*}
where we used the generalized Vandermonde identity $\sum_{\bsm\leq\bsnu,\:|\bsm|=\ell}\binom{\bsnu}{\bsm}=\binom{|\bsnu|}{\ell}$ and $C_P>0$ is the Poincar\'e constant of the embedding $H_0^1(D)\hookrightarrow L^2(D)$. By applying Theorem~\ref{thm:main2} to $G_{\rm nl}(u)$, it follows that the term~\eqref{eq:nonlinearex} decays according to $\mathcal O(s^{-\frac{2}{p}+1})$.
\end{example}

\section{Numerical experiments}\label{sec:numex}
We consider the PDE problem
\begin{align}
\begin{cases}
-\nabla\cdot(a(\bsx,\bsy)\nabla u(\bsx,\bsy))=f(\bsx),&\bsx\in D,~\bsy\in U,\\
u(\cdot,\bsy)=0,&\bsx\in\partial D,~\bsy\in U,
\end{cases}\label{eq:numpde}
\end{align}
over the spatial domain $D=(0,1)^2$ with the source term $f(\bsx)=x_2$. The PDE~\eqref{eq:numpde} is discretized spatially using a finite element method with piecewise linear basis functions and mesh size $h=2^{-5}$.
\subsection{Lognormal input random field}
We let $U=\mathbb R^{\mathbb N}$ and endow the PDE problem~\eqref{eq:numpde} with the lognormal diffusion coefficient
$$
a(\bsx,\bsy)=\exp\bigg(\sum_{j\geq 1}y_j \psi_j(\bsx)\bigg),\quad\bsx\in D,~\bsy\in \mathbb R^{\mathbb N},~\vartheta>1,
$$
where $\psi_j(\bsx) := j^{-\vartheta}\sin(j\pi x_1)\sin(j\pi x_2)$. In this case $(\|\psi\|_{L^{\infty}(D)}) \in \ell^p(\bbN)$ for any $p>\frac{1}{\vartheta}$.
\begin{figure}[!t]
\begin{center}
\includegraphics[width=.7\textwidth]{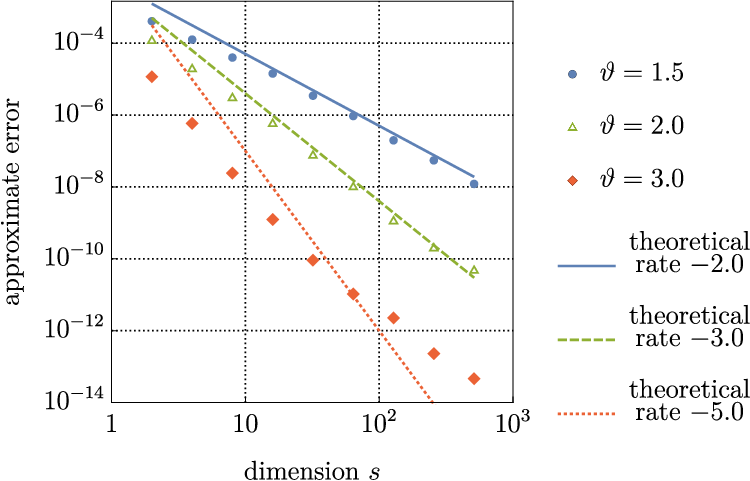}
\end{center}
\caption{The dimension truncation errors corresponding to a lognormally parameterized input random field with decay parameters $\vartheta\in\{1.5,2.0,3.0\}$. The expected dimension truncation error rates are $-2.0$, $-3.0$, and $-5.0$, respectively.}\label{fig:fig}
\end{figure}
To estimate the dimension truncation error, we compute the quantity
\begin{align}
&\bigg\|\int_{\mathbb R^{s'}}(u_{s'}(\cdot,\bsy)-u_s(\cdot,\bsy))\prod_{j=1}^s \varphi_2(y_j)\,{\rm d}\bsy\bigg\|_{H_0^1(D)}\notag\\
&=\bigg\|\int_{(0,1)^{s'}}(u_{s'}(\cdot,\Phi^{-1}(\bst))-u_s(\cdot,\Phi^{-1}(\bst)))\,{\rm d}\bst\bigg\|_{H_0^1(D)},\label{eq:intappx}
\end{align}
where $s'\gg s$ and $\Phi^{-1}$ is the inverse cumulative distribution function of $\prod_{j=1}^{s'}\varphi_2(y_j)$. The high-dimensional integral appearing in~\eqref{eq:intappx} was approximated by using a randomly shifted rank-1 lattice rule with $2^{20}$ cubature nodes and a single random shift. The integration lattice was tailored for each value of the decay parameter $\vartheta$ by using the QMC4PDE software~\cite{qmc4pde,kuonuyenssurvey} and the same random shift was used for each $\vartheta$. As the reference, we use the solution corresponding to $s'=2^{11}$. The value of the high-dimensional integral appearing in~\eqref{eq:intappx} was computed for each finite element node. Since the integrated quantity obtained in this way can also be viewed as an element of the finite element space, it is straightforward to compute its $H_0^1(D)$ norm using the function values at the finite element nodes.

The numerical results are displayed in Figure~\ref{fig:fig} for dimensions $s\in\{2^k: k\in\{1,\ldots,9\}\}$ and decay rates $\vartheta\in\{1.5,2.0,3.0\}$. The corresponding theoretical convergence rates are $-2.0$, $-3.0$, and $-5.0$, respectively, and they are displayed alongside the numerical results. The observed dimension truncation rates corresponding to $\vartheta\in\{1.5,2.0\}$ start off with slower convergence rates and reach the theoretically predicted convergence rates approximately when $s>16$. This behavior appears to be the most extreme in the experiment with $\vartheta=3.0$, where the initial convergence rate for small values of $s$ is slightly slower than the theoretically predicted rate. Moreover, we note that the dimension truncation convergence rate appears to degenerate for large values of $s$ in the case $\vartheta=3.0$, which may be attributed to cubature error when approximating the high-dimensional integral in~\eqref{eq:intappx}.
\subsection{Nonlinear quantity of interest}
Let us revisit Example~\ref{ex:nonlin}. We let $U=[-1,1]^{\mathbb N}$ and endow the PDE problem~\eqref{eq:numpde} with the affine random coefficient
\begin{figure}[!t]
\begin{center}
\includegraphics[width=.7\textwidth]{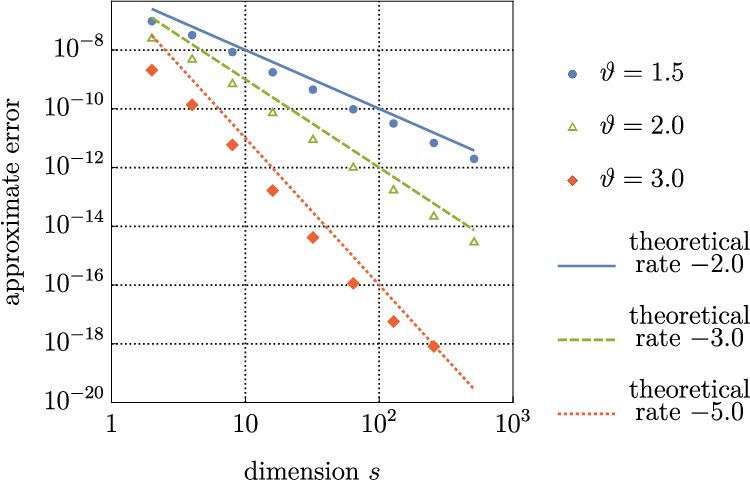}
\end{center}
\caption{The dimension truncation errors corresponding to a nonlinear quantity of interest with decay parameters $\vartheta\in\{1.5,2.0,3.0\}$. The expected dimension truncation error rates are $-2.0$, $-3.0$, and $-5.0$, respectively.}\label{fig:fig2}
\end{figure}
$$
a(\bsx,\bsy)=\frac32+\sum_{j\geq 1}y_j\psi_j(\bsx),\quad\bsx\in D,~\bsy\in [-1,1]^{\mathbb N},~\vartheta>1,
$$
where $\psi_j(\bsx) := j^{-\vartheta}\sin(j\pi x_1)\sin(j\pi x_2)$. In this case $(\|\psi\|_{L^{\infty}(D)}) \in \ell^p(\bbN)$ for any $p>\frac{1}{\vartheta}$.
To estimate the dimension truncation error, we compute
\begin{align}
&\bigg|\int_{[-1,1]^{s'}}(G_{\rm nl}(u_{s'}(\cdot,\bsy))-G_{\rm nl}(u_s(\cdot,\bsy)))2^{-s'}\,{\rm d}\bsy\bigg|\notag\\
&=\bigg|\int_{[0,1]^{s'}}(G_{\rm nl}(u_{s'}(\cdot,2\bsy-1))-G_{\rm nl}(u_s(\cdot,2\bsy-1)))\,{\rm d}\bsy\bigg|,\label{eq:numexint}
\end{align}
where $s'\gg s$ and $G_{\rm nl}$ is the nonlinear quantity of interest defined by~\eqref{eq:nonlinearqoi}. The high-dimensional integrals~\eqref{eq:numexint} were approximated using tailored randomly shifted rank-1 lattice rules with $2^{20}$ cubature nodes and a single random shift, with the same random shift used for each $\vartheta$. The solution corresponding to $s'=2^{11}$ was used as the reference.

The numerical results are displayed in Figure~\ref{fig:fig2} for dimensions $s\in\{2^k: k\in\{1,\ldots,9\}\}$ and decay rates $\vartheta\in\{1.5,2.0,3.0\}$ alongside their respective theoretical convergence rates $-2.0$, $-3.0$, and $-5.0$. The obtained results agree nicely with the theory. Note that the final data point corresponding to the case $\vartheta=3.0$ has been omitted since its difference from the reference solution is smaller than machine precision. %

\section*{Conclusions}
We considered the dimension truncation error analysis for a class of high-dimensional integration problems. The Taylor series approach taken in this paper utilizes the parametric regularity of an integration problem rather than its parametric structure. This should be contrasted with the popular Neumann series approach used to derive dimension truncation error rates for PDEs with random coefficients in the literature. While the utility of the Neumann approach is practically constrained to the so-called affine parametric operator equation setting, our proposed approach appears to be quite robust: with the Taylor series technique, it is possible to prove improved dimension truncation convergence rates even in the non-affine parametric operator equation setting. To this end, we were able to show improved dimension truncation error bounds for elliptic PDE problems with lognormal random diffusion coefficients compared to the standard bounds in the existing literature. Our method also enables the development of dimension truncation rates for smooth nonlinear quantities of interest of the PDE response, provided that the composition of the nonlinear quantity of interest with the PDE solution satisfies the assumptions of our dimension truncation result. Furthermore, since our method has been developed for general separable Banach spaces, our dimension truncation error rates immediately apply for spatially discretized PDE solutions obtained using a conforming finite element method---a property neglected in virtually all of the existing literature on the topic.

\section*{Acknowledgements}
The authors would like to thank Frances Y. Kuo, Claudia Schillings, and Ian H. Sloan for discussions and comments related to this work which helped to improve this paper. Philipp~A. Guth is grateful to the DFG RTG1953 ``Statistical Modeling of Complex Systems and Processes'' for funding of this research. This research includes computations using the computational cluster Katana supported by Research Technology Services at UNSW Sydney~\cite{katana}. The authors also thank the anonymous referees for comments which helped to improve this paper.

\bibliographystyle{plain}
\bibliography{dimtrunc}

\end{document}